\documentclass[11pt]{article}
\usepackage{color}
\usepackage{tikz}
\usepackage{amssymb, amsmath}
\usepackage{epsfig}
\usepackage{graphicx}

\DeclareMathOperator\erf{erf} \DeclareMathOperator\Rea{Re}
\DeclareMathOperator\im{im}

\newtheorem{theorem}{Theorem}

\newtheorem{remark}{Remark}

\newtheorem{definition}{Definition}
\newtheorem{example}{Example}
\newcommand{\proof}{\noindent\textbf{Proof.~}}
\newcommand{\rr}{\mathbb{R}}

\newcommand{\qed}{\space\hfill\hspace*{\fill} $\vbox{\hrule\hbox{\vrule
height1.3ex\hskip1.3ex\vrule}\hrule}$\hss\vskip\topsep\relax}

\begin{document}

\title{A generalization of Laplace and Fourier transforms}

\author{Nikolaos Halidias \\
{\small\textsl{Department of Mathematics }}\\
{\small\textsl{University of the Aegean }}\\
{\small\textsl{Karlovassi  83200  Samos, Greece} }\\
{\small\textsl{email: nikoshalidias@hotmail.com}}}

\maketitle

\begin{abstract} In this note we propose a generalization of the
Laplace and Fourier transforms which we call symmetric Laplace
transform. It combines both the advantages of the Fourier and
Laplace transforms. We give the definition of this generalization,
some examples and basic properties. We also give the form of its
inverse by using the theory of the Fourier transform. Finally, we
apply the symmetric Laplace transform to a parabolic problem and
to an ordinary differential equation.
\end{abstract}

{\bf Keywords} Fourier transform, Laplace transform, functions of
exponential order

{\bf 2010 Mathematics Subject Classification} 42A38, 35A22

\section{Introduction}
In this note we will propose a generalization of the well known
Fourier and Laplace transforms. Let us recall first the notion of
the Fourier transform (see for example \cite{Beerends},
\cite{Dyke}).

\begin{definition} Let $f:\mathbb{R} \to \mathbb{C}$ be a function. We
define a new function $F(s)$ with $s \in \mathbb{R}$ by
\begin{eqnarray*}
F(s) = {\cal F}(f)(s) = \int_{-\infty}^{\infty} f(t) e^{-ist} dt
\end{eqnarray*}
Here $e^{-ix} = \cos x + i \sin x$.
\end{definition}

A common condition in order a function $f$ can be Fourier
transformed is to be absolutely integrable, that is to satisfy the
following definition,
\begin{definition} A function $f: \mathbb{R} \to \mathbb{C}$
is called absolutely integrable on  $\mathbb{R}$ when
\begin{eqnarray*}
\int_{-\infty}^{\infty} | f(t) | dt < \infty
\end{eqnarray*}
The set of all  absolutely integrable functions on $\mathbb{R}$ is
denoted by  $L^1(\rr)$.
\end{definition}
For an absolutely integrable function $f$ the Fourier transform is
well defined since
\begin{eqnarray*}
| F(s) | = \left| \int_{-\infty}^{\infty} f(t) e^{-ist}dt \right|
\leq \int_{-\infty}^{\infty} | f(t) e^{-ist} | dt =
\int_{-\infty}^{\infty} |f(t)| | e^{-ist}| dt
\end{eqnarray*}
But $|e^{-ist}|=1$ so it follows that
\begin{eqnarray*}
|F(s)| \leq \int_{-\infty}^{\infty} |f(t)|dt < \infty
\end{eqnarray*}

A disadvantage of the Fourier transform is that functions like the
Heaviside function, which is the following
\begin{eqnarray*}
H (t) = \left\{%
\begin{array}{ll}
    1, & \mbox{ when } t \geq 0\\[0.2cm]
    0, & \mbox{ when } t < 0\\
\end{array}%
\right.
\end{eqnarray*}
can not be Fourier transformed.

The Laplace transform can be applied to this kind of functions.
Let us recall its definition, (see for example \cite{Schiff}),
\begin{definition}
Let  $f: \rr^{+} \to \mathbb{C}$ and $s \in \mathbb{C}$. We define
the following function
\begin{eqnarray*}
F(s) = {\cal L}(f)(s)  = \int_0^{\infty} e^{-st} f(t)dt
\end{eqnarray*}
where  $s = x+iy$, if the integral exists.
\end{definition}
It is well-known fact that the Laplace transform of the Heaviside
function is the function $\frac{1}{s}$ for $\Rea s > 0$.

A well known condition for a function $f:\rr^+ \to \mathbb{C}$ so
that can be Laplace transformed is to be of exponential order,
that is
\begin{definition} We say that  the function  $f:\rr^{+} \to \mathbb{C}$ is of exponential order  $a$ if there is some constant
 $M > 0$ and some $a>0$ such that for some
 $t_0 \geq 0$, it hold that
\begin{eqnarray*}
|f(t)| \leq M e^{at}, \quad t \geq t_0
\end{eqnarray*}
\end{definition}
A similar definition for some $f:\rr \to \mathbb{C}$ is the
following.
\begin{definition}
We say that the function   $f:\rr \to \mathbb{C}$ is of
exponential order  $a$ if there is some constant $M
> 0$ and some  $a > 0$ such that
\begin{eqnarray*}
|f(t)| \leq M e^{a|t|}, \quad t \in \rr
\end{eqnarray*}
\end{definition}

A disadvantage of the Laplace transform is that it can not be
applied to functions that are defined in all of $\rr$, such as the
following function,
\begin{eqnarray*}
g(t) = \left\{%
\begin{array}{rl}
    1, & \mbox{ when }  t \geq 0\\[0.1cm]
    -1, & \mbox{ when } t < 0 \\
\end{array}%
\right.
\end{eqnarray*}
This function can not be Fourier transformed as well. So, we will
propose a generalization of both previous transforms which we call
the symmetric Laplace transform.

\section{Symmetric Laplace Transform}

We propose the following integral transform,
\begin{definition}
Let $f:\rr \to \mathbb{C}$. We define the symmetric Laplace
transform as follows,
\begin{eqnarray*}
{\cal SL}(f)(x_1,x_2,y) = \int_{-\infty}^{\infty}
e^{(-x_1H(t)+x_2H(-t)-iy)t} f(t) dt
\end{eqnarray*}
where $H(t) = \left\{%
\begin{array}{ll}
    1, & \mbox{ when } t \geq 0\\
    0, & \mbox{ when } t < 0\\
\end{array}%
\right. $
\end{definition}

Choosing $x_1 = x_2 = x$ and $s = x+iy$  we obtain
\begin{eqnarray}\label{simmetrikosLaplace}
{\cal SL}(f(t))(s) = {\cal L}(f(t))(s) + {\cal
L}(f(-t))(\overline{s})
\end{eqnarray}
Obviously, if $f(t) = 0$ for $t < 0$ then the symmetric Laplace
transform coincides with the usual Laplace transform because
${\cal L}(f(-t))(\overline{s}) = 0$. Choosing, $x_1=x_2=0$ then
the symmetric Laplace transform coincides with the Fourier
transform. Therefore, the symmetric Laplace transform generalizes
both the Fourier and the Laplace transforms.

Let us give three examples of the symmetric Laplace transform.

\begin{example}\label{miantistrepsimosLaplace}
We will evaluate the symmetric Laplace transform of the function
\begin{eqnarray*}
f(t) = \left\{%
\begin{array}{rl}
    1, & \mbox{ when }  t \geq 0\\[0.1cm]
    -1, & \mbox{ when } t < 0 \\
\end{array}%
\right.
\end{eqnarray*}
We have that
\begin{eqnarray*}
{\cal SL}(f)(x_1,x_2,y) & = & \int_{-\infty}^{\infty}
 e^{(-x_1H(t)+x_2H(-t)-iy)t} f(t) dt \\ & = & \int_0^{\infty}
 e^{-(x_1+iy)t} dt - \int_{-\infty}^0 e^{-(-x_2+iy)t} dt \\ & = &
 \frac{1}{x_1+iy} + \frac{1}{-x_2+iy}, \quad \mbox{ when } x_1 > 0,
 \; x_2 > 0, \; y \in \rr
\end{eqnarray*}
Choosing  $x_1=x_2=x$ and  $s = x+iy$ then
\begin{eqnarray*}
{\cal SL}(f(t))(s) = \frac{1}{s} - \frac{1}{\overline{s}}, \quad
\mbox{ when  } s = x+iy, \quad x > 0, \; y \in \rr
\end{eqnarray*}
\end{example}

\begin{example} We will compute the symmetric Laplace transform of the function  $f(t) = 1$ with $t \in \rr$.
We have that
\begin{eqnarray*}
{\cal SL}(f)(x_1,x_2,y) & = & \int_{-\infty}^{\infty} e^{(-x_1
H(t) + x_2H(-t) - iy)t} dt \\ & = & \int_{0}^{\infty}
e^{-(x_1+iy)t} dt + \int_0^{\infty} e^{-(x_2-iy)t} dt \\ & = &
\frac{1}{x_1+iy} + \frac{1}{x_2-iy}, \quad \mbox{ for } x_1>0,\;
x_2
> 0, \; y \in \rr
\end{eqnarray*}
Choosing  $x_1=x_2=x$ and $s = x+iy$ it follows that
\begin{eqnarray*}
{\cal SL}(f)(s) = \frac{1}{s} + \frac{1}{\overline{s}}, \quad x >
0, \; y \in \rr
\end{eqnarray*}
\end{example}

\begin{example}
We will compute the symmetric Laplace transform of the following
function
\begin{eqnarray*}
f(t) = \left\{%
\begin{array}{ll}
    \sin xt, & \mbox{ when }  t \geq 0\\[0.2cm]
    \cos xt, & \mbox{ when }  t < 0\\
\end{array}%
\right.
\end{eqnarray*}
We have that
\begin{eqnarray*}
{\cal SL}(f)(x_1,x_2,y) & = & \int_{-\infty}^{\infty}
e^{(-x_1H(t)+x_2H(-t)-iy)t} f(t) dt \\ & = & \int_0^{\infty}
e^{-(x_1+iy)t} \sin xt dt + \int_{-\infty}^0 e^{(x_2-iy)t} \cos xt
dt \\ & = & \frac{x}{(x_1+iy)^2+x^2} +
\frac{x_2-iy}{(x_2-iy)^2+x^2}
\end{eqnarray*}
for $x_1 > 0, \; x_2
> 0, \; y \in \rr$. Choosing  $x_1 = x_2 = z$ and $s = z + iy$
we obtain that
\begin{eqnarray*}
{\cal SL}(f)(s) = \frac{x}{s^2+x^2} +
\frac{\overline{s}}{\overline{s}^2 + x^2}, \quad \mbox{ when }
\Rea s= z
> 0, \; y \in \rr
\end{eqnarray*}

We will also compute the symmetric Laplace transform of the
function
\begin{eqnarray*}
f(t) = \left\{%
\begin{array}{ll}
    \cos xt, & \mbox{ when  }  t \geq 0\\[0.2cm]
    \sin xt, & \mbox{ when }  t < 0\\
\end{array}%
\right.
\end{eqnarray*}
Choosing $x_1=x_2=z$ and $s = z+ iy$ we obtain
\begin{eqnarray*}
{\cal SL}(f)(x,y) & = &  {\cal L}(\cos xt)(s) + {\cal L}(\sin
(-xt)) (\overline{s}) \\ & = & \frac{s}{s^2+x^2} -
\frac{x}{\overline{s}^2+x^2}, \quad \Rea s = z > 0, \; y \in \rr
\end{eqnarray*}
\end{example}

\section{The Inverse of the Symmetric Laplace Transform}
In order to prove that the symmetric Laplace transform is
invertible we will use the following theorem, (see theorem  19.3,
page 248, \cite{Bak}).

\begin{theorem}\label{monadikotitametasximatismou} If the function
$f$ is absolutely integrable on $\rr$ and such that
\begin{eqnarray*}
\int_{-\infty}^{\infty} f(t) e^{-ist} dt = 0, \quad \mbox{ for
every  } s \in \rr
\end{eqnarray*}
then $f(t) = 0$ for every $t \in \rr$.
\end{theorem}

Let us define the set of all the continuous functions $f$ which
are also such that the function $$g(t) = e^{-(x_1H(t) +
x_2H(-t))t}
 f(t)$$
is absolutely integrable for $(x_1,x_2) \in U \subseteq
 \rr^2$ with $U$ chosen appropriately for the given function $f$. We denote this set by  $C_U{\cal SL}$.
For these functions the symmetric Laplace transform is well
defined because ${\cal SL}(f)(s) = {\cal F}(g)(s)$ for $f \in
C_U{\cal SL}$.

The symmetric Laplace transform
 \begin{eqnarray*}
 {\cal SL} : C_U{\cal SL} \to \im {\cal SL}
 \end{eqnarray*}
is a linear transform, recalling the linearity of the integral.
The kernel of this transform contains only the zero element.
Indeed, if
\begin{eqnarray*}
\int_{-\infty}^{\infty} g(t) e^{-ist} dt = 0, \quad \mbox{ for
every } s \in \rr
\end{eqnarray*}
then $g(t) = 0$ for all $ t \in \rr$ and therefore $f(t) = 0$ for
all $t \in \rr$. That means that the symmetric Laplace transform
is a  1-1 linear transform and so invertible. The inverse of this
transform is also linear. These useful conclusions comes from
basic linear algebra (see for example \cite{Berberian},
\cite{Lang} and \cite{Larson}).

However, if we choose $s \in \rr^{+}$ (and not $s \in \mathbb{C}$)
in the symmetric Laplace transform then it follows that is not
invertible. To see this let us compute the symmetric Laplace
transform (with $s \in \rr^{+}$) of the continuous function $f(t)
= t$ with $t \in \rr$. We have
\begin{eqnarray*}
{\cal SL}(f)(s) = \int_0^{\infty} t e^{-st} dt + \int_{-\infty}^0
t e^{st}dt, \quad s > 0
\end{eqnarray*}
It follows that ${\cal SL}(f)(s) = 0$ and thus the kernel of the
transform contains at least one non zero continuous function. That
means that is not invertible.  Therefore, it is important to
choose $s \in \mathbb{C}$ and not in $\rr^{+}$.

In the next theorem we will see the actual form of the inverse of
the symmetric Laplace transform if we assume further that is $f$
is piecewise smooth.

 \begin{theorem}\label{morfiantistrofousimmetrikouLaplace} Let  $f:\rr \to \mathbb{C}$ a piecewise smooth function
 and such that  the function  $g(t) = e^{-(x_1H(t) + x_2H(-t))t}
 f(t)$ is absolutely integrable on $\rr$ when
  $(x_1,x_2) \in U \subseteq
 \rr^2$. If  $F(x_1,x_2,y)$ is the symmetric Laplace transform of $f$ then
 \begin{eqnarray*}
 \lim_{A \to \infty} \frac{1}{2 \pi} \int_{-A}^A F(x_1,x_2,y) e^{ (x_1H(t) - x_2H(-t) + iy)t} dy =
 \frac{1}{2} \left( f(t+) + f(t-) \right), \; (x_1,x_2) \in U
 \end{eqnarray*}
 for every  $t \in \rr$.
 \end{theorem}
\proof We will use the form of the inverse Fourier transform of
the function $g$ (see Theorem 7.3, page 169 of \cite{Beerends}).

Since the function $g(t)$ is absolutely integrable on $\rr$ then
the Fourier transform is well defined and moreover
 ${\cal F}(g)(y) = {\cal
LS}(f)(x_1,x_2,y)= F(x_1,x_2,y)$. Therefore, knowing the form of
the inverse of the Fourier transform we obtain
\begin{eqnarray*}
\frac{1}{2 \pi} \lim_{A \to \infty} \int_{-A}^A F(x_1,x_2,y)
e^{iyt} dy = \frac{1}{2} (g(t+) + g(t-))
\end{eqnarray*}
For  $t > 0$ we have that  $g(t+) = f(t+) e^{-x_1 t}$ and $g(t-) =
f(t-) e^{-x_1t}$ while for  $t < 0$ it holds that  $g(t+) = f(t+)
e^{x_2 t}$ and $g(t-) = f(t-)e^{x_2 t}$. For $t = 0$ we have
$g(0+) = f(0+)$ and $g(0-) = f(0-)$ thus we can write  $g(t+) =
f(t+) e^{(-x_1 H(t) + x_2H(-t))t}$ and $g(t-) = f(t-) e^{(-x_1
H(t) + x_2H(-t))t}$ for $t \in \rr$. That is, it holds that
\begin{eqnarray*}
\frac{1}{2 \pi} \lim_{A \to \infty} \int_{-A}^A
F(x_1,x_2,y)e^{(x_1 H(t) - x_2H(-t)+iy)t}dy = \frac{1}{2} \left(
f(t+) + f(t-) \right), \; (x_1,x_2) \in U
\end{eqnarray*}
for every $t \in \rr$. \qed

\begin{remark}
As we can see the symmetric Laplace transform of a function $f$,
when we choose $x_1=x_2$
 (see
\ref{simmetrikosLaplace}), is the sum of two functions,
$g_1(s)+g_2(\overline{s})$, with $s = x+iy$. The function $g_1(s)$
is the Laplace transform of   $f(t)$ for $t
> 0$ while the function  $g_2(\overline{s})$ is the Laplace transform of  $f(-t)$ for $t > 0$.
Therefore, if we want to find the inverse symmetric Laplace
transform of a function, we should separate it in a sum of two
functions. The first will contain the terms with $s$ and the
second the terms with $\overline{s}$, say $g_1(s)$ and
$g_2(\overline{s})$. Next, we find the inverse of the function
$g_1(s)$ which is equal to $f$ for $t \geq 0$ while the inverse of
the function $g_2(\overline{s})$ is the function $f(t)$ for $t <
0$. Consequently, we have that
\begin{eqnarray*}
f(t) = \left\{%
\begin{array}{ll}
    {\cal L}^{-1}(g_1(s))(t), & \mbox{ when } t \geq 0 \\[0.3cm]
    {\cal L}^{-1}(g_2(\overline{s}))(-t), & \mbox{ when } t < 0 \\
\end{array}%
\right.
\end{eqnarray*}
\end{remark}

\begin{example}
We will find the inverse symmetric Laplace transform of the
function
\begin{eqnarray*}
\frac{1}{s^2} - \frac{1}{\overline{s}^2}
\end{eqnarray*}
Here $g_1(s) = \frac{1}{s^2}$ and $g_2(\overline{s}) =
-\frac{1}{\overline{s}^2}$. Therefore the function $f(t)$ is such
that
\begin{eqnarray*}
f(t) = \left\{%
\begin{array}{ll}
    {\cal L}^{-1}(g_1(s))(t), & \mbox{ when } t \geq 0 \\[0.3cm]
    {\cal L}^{-1}(g_2(\overline{s}))(-t), & \mbox{ when } t < 0 \\
\end{array}%
\right.
\end{eqnarray*}
But, ${\cal L}^{-1}(g_1(s))(t) = t$ and ${\cal
L}^{-1}(g_2(\overline{s}))(-t) =t$ therefore $f(t) = t$ with $t\in
\rr$.
\end{example}

\section{Basic Properties of the Symmetric Laplace Transform}
We will give some basic properties concerning the symmetric
Laplace transform of the derivative of a function.

\begin{theorem} Let  $f ,f^{'}: \rr \to \mathbb{C}$ are continuous  $\rr$
and that  $f$ is of exponential order  $a$. Then
\begin{eqnarray*}
{\cal SL}(f^{'})(s) = s {\cal L}(f(t))(s) - \overline{s} {\cal
L}(f(-t)) (\overline{s})
\end{eqnarray*}
where  $s=x+iy$ and $\Rea s = x > a$, $y \in \rr$.
\end{theorem}
\proof We have that
\begin{eqnarray*}
{\cal SL}(f^{'})(s) = \int_{0}^{\infty} f^{'}(t) e^{-st} dt +
\int_{-\infty}^0 f^{'}(t) e^{\overline{s}t} dt
\end{eqnarray*}
But
\begin{eqnarray*}
\int_{-\infty}^0 f^{'}(t) e^{\overline{s}t} dt & = & \left[ f(t)
e^{\overline{s}t}\right]_{-\infty}^0 - \overline{s}
\int_{-\infty}^0 f(t) e^{\overline{s}t} dt \\ & = &f(0-)-
\overline{s} {\cal L}(f(-t))(\overline{s})
\end{eqnarray*}
and similarly
\begin{eqnarray*}
\int_0^{\infty} f^{'}(t) e^{-st} dt & = & \left[ f(t)
e^{-st}\right]_{0}^{\infty} + s \int_{0}^{\infty} f(t) e^{-st} dt
\\ & = & -f(0+)+ s {\cal L}(f(t))(s)
\end{eqnarray*}
Using the continuity of  $f$ we get the desired result.
 \qed
Similarly, we have the following result.
\begin{theorem}
 Suppose that  $f,f^{'},f^{''}$ are continuous on $\rr$ and that
 $f,f^{'}$ are of exponential order  $a$. Then
 \begin{eqnarray*}
 {\cal SL}(f^{''})(s) = s^2 {\cal L}(f(t))(s)+\overline{s}^2 {\cal
 L}(f(-t)) (\overline{s}) - f(0)(s+\overline{s})
 \end{eqnarray*}
 where  $s = x+iy$ and $x > a, y\in \rr$.
 \end{theorem}

More general result in this direction is the following using
induction.

\begin{theorem}
Let $f,f^{'},\cdots,f^{(n)}$ are continuous on  $\rr$, except
maybe at zero, while the functions $f,f^{'},\cdots,f^{(n-1)}$ are
of exponential order  $a$. Then
\begin{eqnarray*}
& & {\cal SL}(f^{(n)}(t))(s) \\[0.3cm]  &  =  & {\cal
L}\left(f^{(n)}(t)\right)(s)+ {\cal
L}\left(f^{(n)}(-t)\right)(\overline{s})
\\[0.3cm] &  = & \underbrace{s^n {\cal L}(f) -s^{n-1}f(0+) - s^{n-2}
f^{'}(0+) - \cdots - f^{(n-1)}(0+)}_{{\cal L}(f^{(n)}(t))(s)} \\[0.5cm] &
& + \underbrace{(-\overline{s})^n {\cal L}(f(-t))(\overline{s}) +
(-\overline{s})^{n-1} f(0-) + (-\overline{s})^{n-2} f'(0-) +
\cdots + f^{(n-1)}(0-)}_{{\cal
L}\left(f^{(n)}(-t)\right)(\overline{s})}
\end{eqnarray*}
\end{theorem}

\begin{example}
We will study the following parabolic problem,
\begin{eqnarray*}
\begin{array}{lllll}
u_{xx}(x,t) & = & u_t(x,t), \quad & x \in \rr, \quad &t > 0
\\[0.1cm]
u(x,0) &= & f(x), \quad & x\in \rr \\[0.1cm]
u(0,t) & = & 0, \quad & & t>0 \\[0.1cm]
\end{array}
\end{eqnarray*}
Here  $f$ is as follows
\begin{eqnarray*}
f(x) = \left\{%
\begin{array}{rl}
    1, & \mbox{ when } x \geq 0 \\[0.1cm]
    -1, & \mbox{ when } x < 0 \\
\end{array}%
\right.
\end{eqnarray*}
We will apply the symmetric Laplace transform on the $x$ variable
to get
\begin{eqnarray}\label{metasximatismenisimmetrika}
s^2 G(s,t) + \overline{s}^2 \tilde{G}(\overline{s},t) = G_t(s,t) +
\tilde{G}_t(\overline{s},t)
\end{eqnarray}
where $G(s,t) = {\cal L}(u(x,t))$ and  $\tilde{G}(s,t) = {\cal
L}(u(-x,t))$.

We apply also the symmetric Laplace transform to the condition
$u(x,0) = f(x)$ to get
\begin{eqnarray}\label{simetrikessinthikes}
G(s,0) + \tilde{G}(\overline{s},0) = \frac{1}{s} -
\frac{1}{\overline{s}}
\end{eqnarray}
From  \ref{metasximatismenisimmetrika} and
\ref{simetrikessinthikes} we get two ordinary differential
equations which are
\begin{eqnarray*}
\begin{array}{lllrrr}
s^2 G(s,t) & = & G_t(s,t), \quad G(s,0) & = & \frac{1}{s},
\\[0.3cm]
\overline{s}^2 \tilde{G}(\overline{s},t) & = &
\tilde{G}_t(\overline{s},t), \quad \tilde{G}(\overline{s},0) &  =&
 - \frac{1}{\overline{s}}
\end{array}
\end{eqnarray*}
Their solutions are
\begin{eqnarray*}
G(s,t) = \frac{1}{s} e^{s^2t}, \quad \tilde{G}(\overline{s},t) = -
\frac{1}{\overline{s}} e^{\overline{s}^2 t}
\end{eqnarray*}
Inverting the symmetric Laplace transform (using the convolution
theorem of the usual Laplace transform) we arrive at
\begin{eqnarray*}
u(x,t) = \left\{%
\begin{array}{rl}
    \erf\left(\frac{x}{2\sqrt{t}}\right), & \mbox{ when } x \geq 0 \\[0.5cm]
   - \erf\left(\frac{-x}{2\sqrt{t}}\right), & \mbox{ when } x < 0 \\
\end{array}%
\right.
\end{eqnarray*}
where $\displaystyle \erf\left(\frac{x}{2\sqrt{t}}\right) =
\frac{2}{\sqrt{\pi}} \int_0^{\frac{x}{2\sqrt{t}}} e^{-u^2} du$. It
is  easy  to show that the function  $u(x,t)$ satisfies our
problem.

Next, assuming further that every solution of our parabolic
problem should be bounded as
 $|x| \to \infty$ while  $u_x(x,t) \to 0$ as $|x|
\to \infty$ we will prove that this problem admits a unique
solution. Note that these extra conditions satisfied by the
solution above.

Suppose that the problem admits at least two solutions, the
$u_1,u_2$. Then their difference $v = u_1-u_2$ will satisfy the
following problem
\begin{eqnarray*}
\begin{array}{lllll}
v_{xx}(x,t) & = & v_t(x,t), \quad & x \in \rr, \quad &t > 0
\\[0.1cm]
v(x,0) &= & 0, \quad & x\in \rr \\[0.1cm]
v(0,t) & = & 0, \quad & & t>0 \\[0.1cm]
v(x,t) & \mbox{ bounded as } & |x| \to \infty, \quad & v_{x}(x,t)
\to 0 \mbox{ as } & |x| \to \infty
\end{array}
\end{eqnarray*}
Multiply now the equation by
 $v(x,t)$ and integrate over $(0,t)$. On the right hand side of
 the equation we get
\begin{eqnarray*}
\int_0^t v_t(x,y) v(x,y) dy  = v^2(x,t) - \int_0^t   v(x,y)
v_t(x,y) dy
\end{eqnarray*}
Therefore
\begin{eqnarray*}
\int_0^t v_t(x,y) v(x,y) dy = \frac{1}{2} v^2(x,t)
\end{eqnarray*}
Next, we integrate over $\rr$ to get on the left side of the
equation,
\begin{eqnarray*}
\int_{-\infty}^{\infty} \int_0^t v_{xx}(y,r) v(y,r) dr dy & = &
\int_0^t \int_{-\infty}^{\infty}  v_{xx}(y,r) v(y,r) dy dr \\ & =
& \int_0^t \left( v_x(y,r)v(y,r)\Big|_{y=-\infty}^{y=\infty}
\right) - \int_{-\infty}^{\infty} v_x^2(y,r) dydr \\ & = & -
\int_0^t \int_{-\infty}^{\infty} v_x^2(y,r) dydr
\end{eqnarray*}
Therefore, it holds that
\begin{eqnarray*}
\frac{1}{2}\int_{-\infty}^{\infty} v^2(y,t) dy + \int_0^t
\int_{-\infty}^{\infty} v_x^2(y,r) dydr = 0
\end{eqnarray*}
which is true only when $v(x,t) = 0$ for every  $x \in \rr$ and $t
> 0$, that is  $u_1(x,t) = u_2(x,t)$, therefore the problem has a unique solution.
\end{example}

\begin{example}
We will evaluate the solution of the following ordinary
differential equation
\begin{eqnarray*}
y^{''}(t) + y(t) & = & f(t) = \left\{%
\begin{array}{ll}
    e^t, & \mbox{ όταν }  t \geq 0\\[0.2cm]
    1, & \mbox{ όταν } t < 0 \\
\end{array}%
\right.    , \quad t \in \mathbb{R} \\
y(0) & = & 0
\end{eqnarray*}
Obviously, we can not apply the Fourier transform neither the
usual Laplace transform. We can apply the symmetric Laplace
transform to get,
\begin{eqnarray*}
s^2 {\cal L}(y(t))(s) + \overline{s}^2 {\cal
L}(y(-t))(\overline{s}) + {\cal L}(y(t))(s) + {\cal
L}(y(-t))(\overline{s}) = \frac{1}{s-1} + \frac{1}{\overline{s}},
\quad \Rea s > 1
\end{eqnarray*}
Equating similar terms, we get
\begin{eqnarray*}
{\cal L}(y(t))(s) & = & \frac{1}{2} \frac{1}{s-1} - \frac{1}{2}
\frac{s}{s^2+1} - \frac{1}{2} \frac{1}{s^2+1} \\
{\cal L}(y(-t))(\overline{s}) & = & \frac{1}{\overline{s}} -
\frac{\overline{s}}{\overline{s}^2+1}
\end{eqnarray*}
 Computing the inverse transforms we arrive at
 \begin{eqnarray*}
 y(t) = \left\{%
\begin{array}{cl}
    \frac{1}{2} e^t - \frac{1}{2} \cos t - \frac{1}{2} \sin t, & \mbox{ όταν } t \geq 0
    \\[0.3cm]
   1 -\cos t, & \mbox{ όταν } t < 0 \\
\end{array}%
\right.
\end{eqnarray*}
It is easy to see that $y(t)$ satisfy the above ode.
\end{example}

\end{document}